\documentclass[a4paper,11pt]{amsart}

\usepackage[english]{my-shortcuts}
\usepackage{srcltx,algorithm,algorithmic}
\usepackage{a4wide}
\usepackage{url}

\begin{document}

\title[Spacings of Laguerre and Wishart]
{On the spacings between the successive zeros of the Laguerre polynomials}
\author{St\'ephane Chr\'etien {\tiny and} S\'ebastien Darses}

\address{Laboratoire de Math\'ematiques, UMR 6623\\ 
Universit\'e de Franche-Comt\'e, 16 route de Gray\\
25030 Besancon, France} 
\email{stephane.chretien@univ-fcomte.fr}

\address{I2M, UMR 6632\\
Aix-Marseille Universit\'e, Technop\^ole Ch\^{a}teau-Gombert\\
39 rue Joliot Curie\\ 13453 Marseille Cedex 13, France\\
and 
Universit\'e de Franche-Comt\'e, 16 route de Gray\\
25030 Besancon, France}
\email{sebastien.darses@univ-amu.fr}

\maketitle


\begin{abstract}
We propose a simple uniform lower bound on the spacings between the successive zeros of 
the Laguerre polynomials $L_n^{(\alpha)}$ for all $\alpha>-1$. 
Our bound is sharp regarding the order of dependency on $n$ and $\alpha$ in various ranges.
In particular, we recover the orders given in \cite{ahmed} for $\alpha \in (-1,1]$.  
 \end{abstract}


\section{Introduction}

The study of orthogonal polynomials has a long history with exciting interplay with numerous fields, including random matrix theory. The Laguerre polynomials which occur as the solutions of important differential equations \cite{Szego:AMS75}, have had many applications in 
physics (electrostatics, quantum mechanics \cite{FreedenGutting:Birkhauser13}), engineering (control theory; see e.g. \cite{DattaMohan:WS95}), random matrix theory (Wishart distribution; 
see e.g. \cite{DetteImhof:TAMS07} and \cite{Faraut:LectNotes11}) and many other fields. 
The knowledge of the spacings between successive zeros of the Laguerre polynomials, interesting in its own right, 
is also potentially of great interest in many situations, e.g. for the spacings between successive eigenvalues 
of Wishart matrices, for bounding the gaps between sucessive energy levels in quantum mechanics or for the analysis 
of numerical algorithms in system identification problems, to name a few.  

In this short note, we provide a uniform lower bound for the gaps between successive zeros of the Laguerre polynomials 
$L_n^{(\alpha)}$. In \cite{ahmed}, important bounds were proposed in the case $\alpha\in(-1,1]$ for individual spacings 
(i.e. bounds depending also on the ranking). Our bound is uniform but it is valid on the entire range $\alpha>-1$. 
For this reason, our bound might be helpful in a large number of applications. In particular, the cases including large values of $\alpha$, are those of interest for random matrices with Wishart distribution.
Our approach is based on a remarkable well known identity (a Bethe {\em ansatz} equation; see e.g. \cite{Krasikov},\cite{Krasovsky}). 
 
\section{Preliminaries: Bethe ansatz equality}

We first recall the following remarkable general result, see e.g. Lemma 1 in \cite{Krasikov}.
Let $f$ be a polynomial with real simple zeros $x_1<\cdots <x_n$, satisfying the ODE $f''-2af'+bf=0$ where $a$ and $b$ 
are meromorphic function whose poles are different from the $x_i$'s. Then for any fixed $k\in\{1\cdots n\}$,
\bea  \label{bethe}
\sum_{j\neq k} \frac{1}{(x_k-x_j)^2} & = & 
\frac{\Delta(x_k)-2 a'(x_k)}{3},
\eea 
with $\Delta(x)=b(x)-a^2(x)$. Such equalities are called Bethe {\em ansatz} equations.

For $\alpha>-1$, the Laguerre polynomials $L_n^{(\alpha)}$ ($n$ indicates the degree) are 
orthogonal polynomials with respect to the weight $x^\alpha e^{-x}$ on $(0,\infty)$. Let 
$x_{n,n}(\alpha) < \cdots < x_{n,1}(\alpha)$ 
denote the zeros of $L_n^{(\alpha)}$. It is known that the polynomial $L_n^{(\alpha)}$ is a solution of the second order ODE:
\bean
u'' - \left(1-\frac{\alpha+1}{x}\right)u' +\frac{n}{x} u & = & 0.
\eean
In this case, $a(x)=\frac12 \left(1-\frac{\alpha+1}{x}\right)$. Therefore, 
\bean
\Delta(x) =\frac{n}{x}-\frac{(x-\alpha-1)^2}{4x^2}=\frac{-x^2+(2(\alpha+1)+4n)x-(\alpha+1)^2}{4x^2},
\eean
and then using the notations in \cite{Krasikov},
\bea \label{Delta}
\Delta(x) & = & \frac{(U^2-x)(x-V^2)}{4x^2},
\eea
where 
\beq \label{uv}
U =  \sqrt{n+\alpha+1}+\sqrt{n}, \quad V = \sqrt{n+\alpha+1}-\sqrt{n}.
\eeq

Since the l.h.s. of (\ref{bethe}) is positive and $a'(x)>0$ for $x>0$, an immediate consequence of (\ref{bethe})
is that for all $k$, $(U^2-x_{n,k}(\alpha))(x_{n,k}(\alpha)-V^2)>0$, i.e.
\bea \label{extreme}
V^2 < x_{n,n}(\alpha)  < x_{n,1}(\alpha) < U^2.
\eea 
Several bounds for the extreme zeros are known and can be found in \cite{DimitrovNikolov:JAP10,Gatteschi,Ismail,Krasikov,Szego:AMS75}. 
For instance, using the Bethe {\em ansatz}, Krasikov proved \cite[Theorem 1]{Krasikov}:
\beq \label{krasikov}
V^2 +3V^{4/3}(U^2-V^2)^{-1/3} \le x_{n,n}(\alpha)  < x_{n,1}(\alpha) \le U^2 -3U^{4/3}(U^2-V^2)^{-1/3}+2.
\eeq

\section{Main result}

We show by means of elementary computations that the Bethe {\em ansatz} equality actually yields a simple uniform lower bound for
$x_{n,k}(\alpha)-x_{n,k+1}(\alpha)$, which turns out to be sharp, see Remark (2) below.

\begin{theo}\label{main} Let $\alpha>-1$. 
Then, the following lower bound for the spacings holds for all $k\in\{1,\cdots, n-1\}$:
\bea \label{mainbound}
x_{n,k}(\alpha)-x_{n,k+1}(\alpha) & \ge &  \sqrt{3} \frac{\alpha+1}{\sqrt{n(n+\alpha+1)}}.
\eea 
Moreover, if $\alpha \ge n/C$ for some $C>0$, we have
\bea \label{rangen}
x_{n,k}(\alpha)-x_{n,k+1}(\alpha) & \ge & \frac{1}{\sqrt{C+1}} \sqrt{\frac{\alpha}{n}}.
\eea
\end{theo}

\smallskip

\subsection{Proof of Theorem \ref{main}}\

From (\ref{bethe}), (\ref{extreme}) and $a'(x)>0$ for $x>0$, we deduce the following inequality
\bea \label{bethe2} 
\frac{1}{(x_{n,k}(\alpha)-x_{n,k+1}(\alpha))^2} \le \sum_{j\neq k} \frac{1}{(x_{n,k}(\alpha)-x_{n,j}(\alpha))^2} \le \frac13 \sup_{V^2\le x \le U^2}\Delta(x).
\eea
The first inequality above seems to be crude, but is not, see Remark (1) below.

Let us then study the function $\Delta$. The derivative of $\Delta$ on $(0,+\infty)$ reads:
\bean
\Delta'(x) = \frac{(-2x+U^2+V^2)x^2-2x\left(-x^2+(U^2+V^2)x-U^2V^2\right)}{4x^4} = \frac{2U^2V^2-(U^2+V^2)x}{4x^3}.
\eean
Thus, $\Delta$ has a unique maximum on $(0,+\infty)$ that is reached at $x^*=\frac{2U^2V^2}{U^2+V^2}$. We have:
\bean
U^2- x^*& = & \frac{U^4- U^2V^2}{U^2+V^2}=U^2 \frac{U^2-V^2}{U^2+V^2} \\
x^*-V^2 & = & \frac{U^2V^2-V^4}{U^2+V^2}=V^2 \frac{U^2-V^2}{U^2+V^2}.
\eean
Thus, we obtain by plugging into (\ref{Delta}),
\bean
\sup_{V^2\le x \le U^2}\Delta(x)=\Delta(x^*) = \frac{(U^2-V^2)^2}{16\ U^2V^2},
\eean
since one can check that $x^*\in(V^2,U^2)$.
Moreover, from the expressions (\ref{uv}) of $U$ and $V$:  
\bean
U^2- V^2 & = & (U-V)(U+V)=4\sqrt{n}\sqrt{n+\alpha+1}\\
UV & = & \alpha+1.
\eean
Hence, plugging these last equalities in (\ref{bethe2}), we can write
\bean 
\frac{1}{(x_{n,k}(\alpha)-x_{n,k+1}(\alpha))^2} & \le & 
\frac{1}{3}\ \frac{4^2\ n(n+\alpha+1)}{16\ (\alpha+1)^2},
\eean 
and finally
\bean
x_{n,k}(\alpha)-x_{n,k+1}(\alpha) & \ge & \sqrt{3} \frac{\alpha+1}{\sqrt{n(n+\alpha+1)}}.
\eean

Now assume that $n\le C\alpha$. Then $n+\alpha+1\le (C+1)\alpha +1$. Therefore $\sqrt{n+\alpha+1}\le \sqrt{2(C+1)\alpha}$, where we used  $1\le C\alpha\le (C+1)\alpha $. Hence
\bean
x_{n,k}(\alpha)-x_{n,k+1}(\alpha) & \ge & \sqrt{\frac{3}{2(C+1)}} \sqrt{\frac{\alpha}{n}},
\eean
which completes the proof of Theorem \ref{main}.

\subsection{Remarks}

\begin{enumerate}
\item 
Notice that replacing the sum $\sum_{j\neq k} (x_{n,k}(\alpha)-x_{n,j}(\alpha))^{-2}$ by the single term 
$(x_{n,k}(\alpha)-x_{n,k+1}(\alpha))^{-2}$ does not deteriorate {\it a priori} the order of dependency on $n$ and $\alpha$ 
of a uniform bound in $k$ of $x_{n,k}(\alpha)-x_{n,k+1}(\alpha)$. Indeed, let $0<\delta<x_{n,k}(\alpha)-x_{n,k+1}(\alpha)$ 
for all $k$, we have the following simple inequality for any fixed $k$:
\bean
\frac1{(x_{n,k}(\alpha)-x_{n,k+1}(\alpha))^2} \le \sum_{j\neq k} \frac{1}{(x_{n,k}(\alpha)-x_{n,j}(\alpha))^2} \le \sum_{j\neq k} \frac{1}{(\delta|j-k|)^2}
\le 2\frac{\pi^2}{6}\  \frac{1}{\delta^2}.
\eean

\item \label{rem2} Let us verify that our bound is sharp regarding the order of dependency on $n$ and $\alpha$ 
 in various ranges.\\
{\bf Case $\alpha\in(-1,1]$}: Theorem 5.1 in \cite{ahmed} says that for all $\alpha\in(-1,1]$:
\bean
\left(n+(\alpha+1)/2\right) \left(x_{n,k}(\alpha)-x_{n,k+1}(\alpha)\right) & \xrightarrow[]{n\to \infty} & j_{\alpha,k+1}^2-j_{\alpha,k}^2, 
\eean
where $j_{\alpha,k}$ is the $k$-th zeros of the Bessel function $J_\alpha(x)$.
But, for all $k\ge 1$, the following holds (see \cite[Theorem 3]{Hethcote} and \cite[p.2]{Finch}):
\bean
\pi \le j_{\alpha,k+1}-j_{\alpha,k} \le  2\pi  \\
j_{\alpha,k+1}+j_{\alpha,k} \ge  2\sqrt{(k-1/4)^2\pi + \alpha^2} \ge 1+\alpha.
\eean
As a consequence, for small $k$, $x_{n,k}(\alpha)-x_{n,k+1}(\alpha) \sim C(\alpha)/n$, which is consistent with our bound (\ref{mainbound}).\\
{\bf Case $n\le C\alpha$} for an absolute constant $C>0$: Summing (\ref{rangen}) over $k$ yields
\bean
\frac{\sqrt{n\alpha}}{\sqrt{C+1}}  & \le &\sum_{1\le k \le n-1} \left(x_{n,k}(\alpha)-x_{n,k+1}(\alpha)\right) = x_{n,1}(\alpha)-x_{n,n}(\alpha) \\
& \le & U^2- V^2 =4\sqrt{n}\sqrt{n+\alpha+1} \le 6\sqrt{C+1}\sqrt{n\alpha},\nonumber
\eean
which means that the bound (\ref{rangen}) is sharp with respect to 
the orders of $n$ and $\alpha$ up to a multiplicative constant.

Notice moreover that in full generality, $C$ can be taken as a function of $n$ with absolutely no change in the proof.

\item Finally, since the Bethe {\em ansatz} equation (\ref{bethe}) is a general equality for polynomials $f$ with real simple zeros, satisfying the ODE $f''-2af'+bf=0$, good prior bounds on the 
extreme zeros for such polynomials could be used to obtain similar results as Theorem \ref{main}.

\end{enumerate}

\newpage

\section{Numerical results}

We now provide numerical results on the successive spacings of the Laguerre polynomials $L_n^{(\alpha)}$ for various values of $n$ and $\alpha$. 

\begin{figure}[h!]
\begin{center}
\includegraphics[width=14cm]{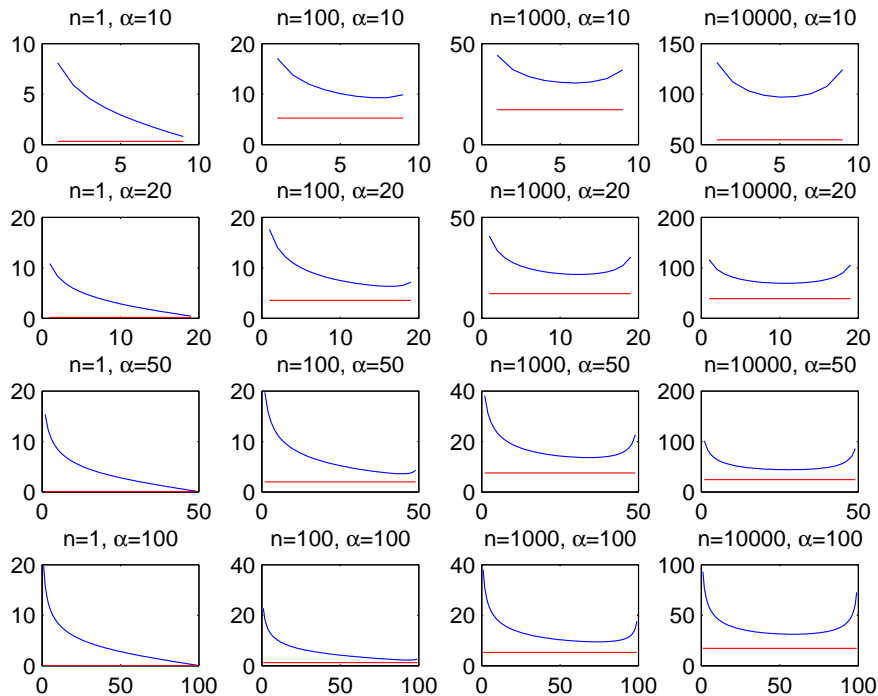}
\caption{Comparison between the uniform bound (\ref{mainbound}) in red, and the function $i\mapsto x_{n,i}(\alpha)-x_{n,i+1}(\alpha)$, $1\le i\le n-1$, in blue. We set $\alpha=1,100,10^3,10^4$ and $n=10,20,50,100$.}\label{gaga}
\end{center}
\end{figure}

Let us make a few comments on Figure \ref{gaga}. The first 
column illustrates that the uniform bound almost coincides with the smallest spacing, which is here $x_{n,n-1}(1)-x_{n,n}(1)$ 
(Recall that $x_{n,n}(\alpha)$ is the smallest zero). When $\alpha$ is large compared to $n$, the behavior 
is quite different. For instance, based on Remark \ref{rem2} in the case $\alpha\ge n/C$,
we can expect most spacings to be almost equal, 
i.e. close to the uniform lower bound $\sqrt{\alpha/n}$ up to a multiplicative constant.
In the last two columns of Figure \ref{gaga} (large values of $\alpha$ compared to $n$), we observe that this phenomena actually occurs in the bulk, i.e. for $\e n\le i \le (1-\e)n$, $0<\e<1$. 

The results plotted in Figure \ref{gaga} have been obtained using {\sc Matlab} and the codes 
available at \url{http://people.sc.fsu.edu/~jburkardt/m_src/laguerre_polynomial/}.

\bibliographystyle{amsplain}
\bibliography{database}

\end{document}